\documentclass[a4paper,11pt]{article}

\usepackage[utf8]{inputenc}
\usepackage[UKenglish]{babel}
\usepackage{mathtools,amssymb,amsthm}
\usepackage{adjustbox,graphicx,float}
\usepackage{hyperref}
\usepackage{listings}
\usepackage{cases}
\usepackage{lineno}
\usepackage{a4wide}
\include{linenopatch}

\newtheorem{prop}{Proposition} 
\theoremstyle{definition}

\newtheorem{thm}[prop]{Theorem}

\newtheorem{defn}[prop]{Definition}

\usepackage{csquotes}
\title{Maximising the number of regions when embedding a $N$-cycle graph.}
\author{Adam Dunajski \\ University of Edinburgh.\\
Email: {\tt s1901024@ed.ac.uk} }
\begin{document}
\maketitle
\begin{abstract}
We find the maximal number of regions that a straight line embedding of a $N$-cycle graph can enclose.
\end{abstract}
\section*{The Problem.}
In this paper we find the maximal number of regions, $f(N)$, that a straight line embedding of a $N$-cycle graph can enclose. This was a problem presented to me by my former teacher Dave Strong, who had also worked out the first few cases by hand. He pointed out that as far as he knew a general answer was not known. I am grateful to Prof. Richard Schwartz and Prof. Andrew Thomason for their guidance and insight into the problem.

\quad

\begin{figure}[H]
\centering
\includegraphics[width=0.5\textwidth]{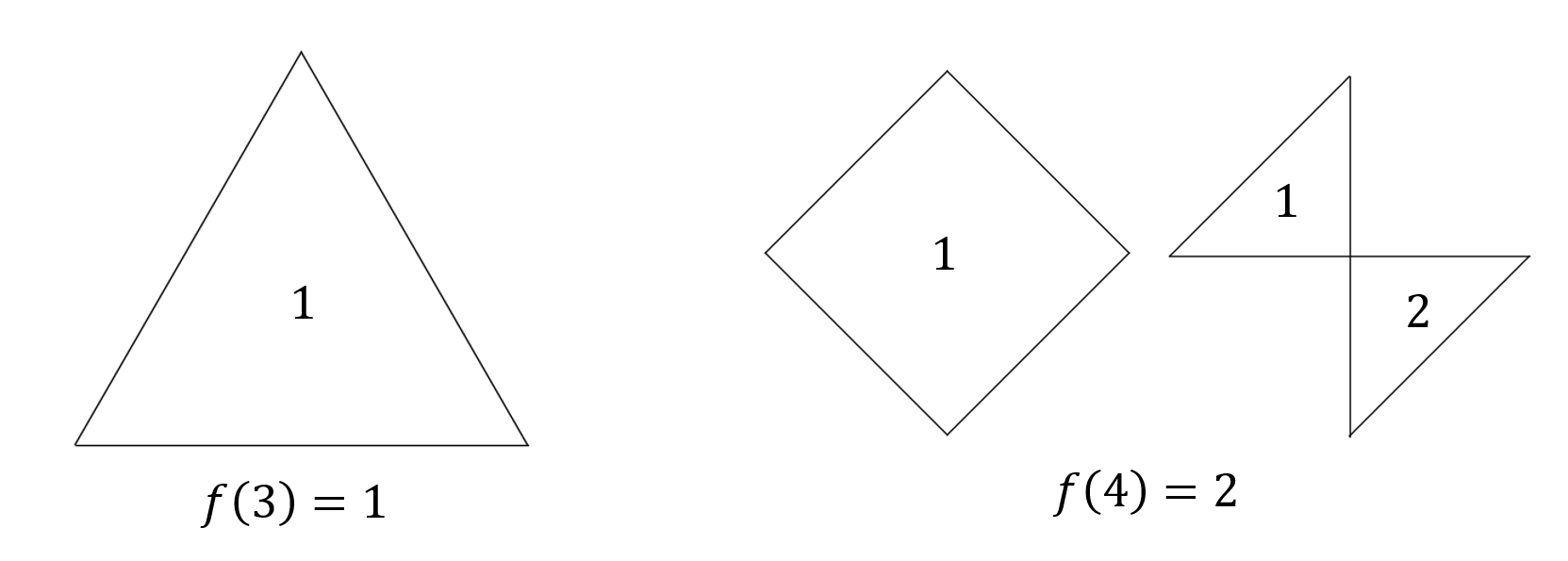}
\caption{$f(3)$ and $f(4)$.}
\label{3_4}
\end{figure}

Recall (see e.g. \cite{Diestel}) that a $N$-cycle graph is a graph that consists of $N$ corners\footnote{\label{note1}We use the term corners instead of vertices. Similarly we use the term line-segment instead of edge for the connections between corners, as the line-segments of a $N$-cycle graph will later be split into edges. For example, the optimal arrangement for $N = 4$ shown in Figure \ref{3_4} we have four corners but five vertices, and four line-segments but six edges.}(at least three, if the graph is simple) connected in a closed chain.  For example a $3$-cycle graph consists of corners $1,2,$ and $3$ connected $1-2-3-1$. When embedding a $3$-cycle graph, one is forced to arrange the corners in some triangle. This encloses one region as shown on the left of Figure \ref{3_4}. We therefore get $f(3) = 1$. When we embed a $4$-cycle graph however we have two options; We could either create a quadrilateral that encloses one region, or have two opposite line-segments cross to create an hourglass shape that encloses two regions. Both of these arrangements are shown on the right of Figure \ref{3_4}. The latter of these encloses a larger number of regions, and as these are the only two possibilities we get $f(4) = 2$.

The objective is then to find a formula, $f(N)$, that gives the maximal number of regions for a $N$-cycle graph. Note that in Figure \ref{3_4} the corners of the cycle graphs are arranged in the shapes of a regular $3$-gon and $4$-gons. This is not a requirement, and in fact for larger $N$ restricting to a regular $N$-gon can lead to triple intersections which reduce the number of regions.

\section*{The Solution.}
We will prove the following.
\begin{thm}
The maximal number of regions created by a straight line embedding of a $N$-cycle graph is given by $f(N)$ defined as follows:
\begin{subnumcases}{f(N)=}
   \frac{1}{2} N^2 - 2 N + 2, & \text{for even} $N$. \label{eventhm}\\
    \frac{1}{2} N^2 - \frac{3}{2} N + 1, & \text{for odd}  $N$. \label{oddthm}
\end{subnumcases}
\label{F}
\end{thm}
Finding $f(N)$ becomes equivalent to maximising the number of intersections between the line-segments in the $N$-cycle. This is a result of Euler's formula relating the number of vertices, edges, and faces. In our case, as we do not include the infinite outer region in the count, we use the following
\begin{equation}
    F = E - V + 1 .
\label{Euler}
\end{equation}
By having two line-segments intersect we increase the number of edges by two, whilst increasing the number of vertices by only one. By (\ref{Euler}) we then see that each additional intersection results in an additional region.

A given line-segment in a $N$-cycle graph intersects any other line-segment at most once. A special role will be played by line-segments with the maximal possible number of intersections. 
\begin{defn}
A \textit{splitter} is a line-segment in a $N$-cycle graph that intersects all other line-segments. Note that this could be at a corner.
\end{defn}
A splitter therefore contains $N-1$ vertices and $N-2$ edges. We say a line-segment is one-off being a splitter if it intersects all but one other line-segment. In this case such a line-segment would have $N-2$ vertices and $N-3$ edges. We refer to these as \textit{one-off splitters}.

We now split the proof of Theorem \ref{F} into two cases, depending on the parity of $N$.

\subsection*{$N$ is odd.}
\begin{figure}[H]
\centering
\includegraphics[width=1\textwidth]{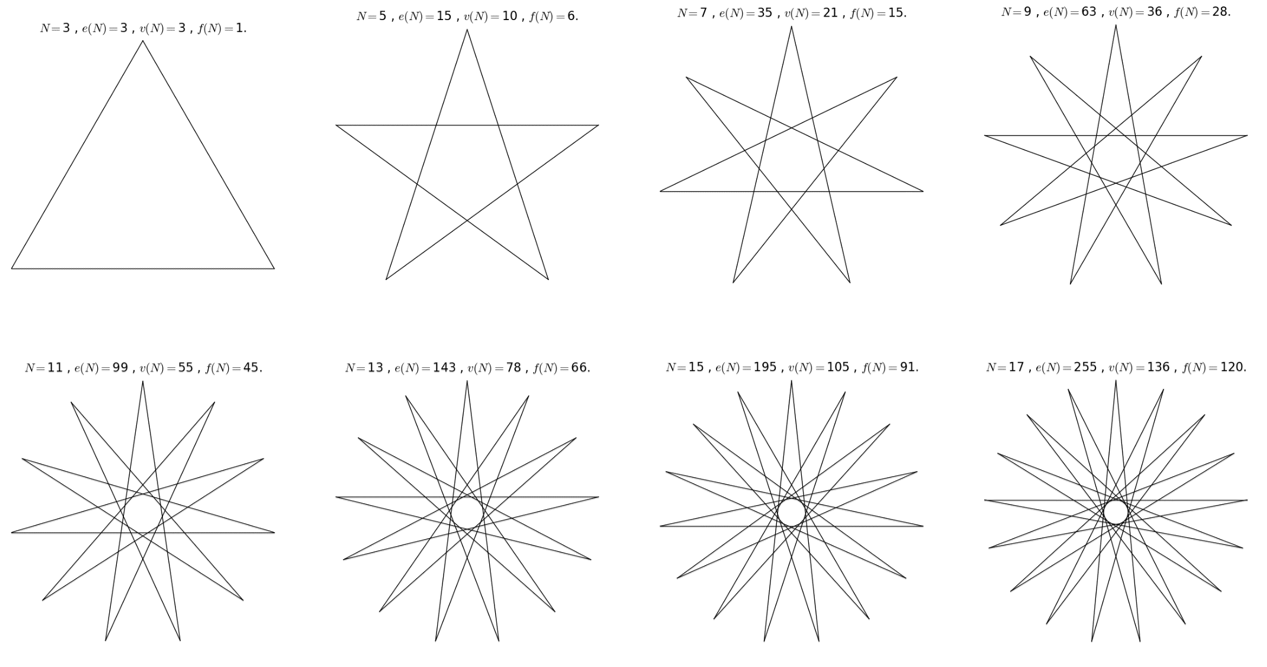}
\caption{Optimal arrangements for odd $N$.}
\label{odd}
\end{figure}

For odd $N$ it is possible for every line-segment to be a splitter. The resulting arrangements are therefore optimal as they have the maximal possible number of vertices. To construct these we label the $N$ corners $1$ to $N$. We are inspired by the pentagram, where $N = 5$, and connect corner $c$ with corner $\big(c + \frac{N-1}{2}\big)$ mod $N$. This results in the arrangements shown in Figure \ref{odd}. 
We then work out the number of regions as follows. Let $e(N)$ be the number of edges in the optimal arrangement and $v(N)$ be the number of vertices. We know that each of the $N$ line-segments contains $N-1$ vertices. However each vertex belongs to exactly two line-segments. This gives us
\[
v(N) = \frac{N(N-1)}{2}.
\]
Similarly each of the $N$ line-segments contains $N-2$ edges. This gives us 
\[
e(N) = N(N-2).
\]
We then combine these using Euler's formula i.e. (\ref{Euler}) to give us 
\begin{align*}
  f(N) &= e(N) - v(N) + 1 \\
  &=  \frac{1}{2} N^2 - \frac{3}{2} N + 1.
\end{align*}
This concludes the proof of (\ref{oddthm}).

\subsection*{ $N$ is even.}
\begin{figure}[H]
\centering
\includegraphics[width=0.5\textwidth]{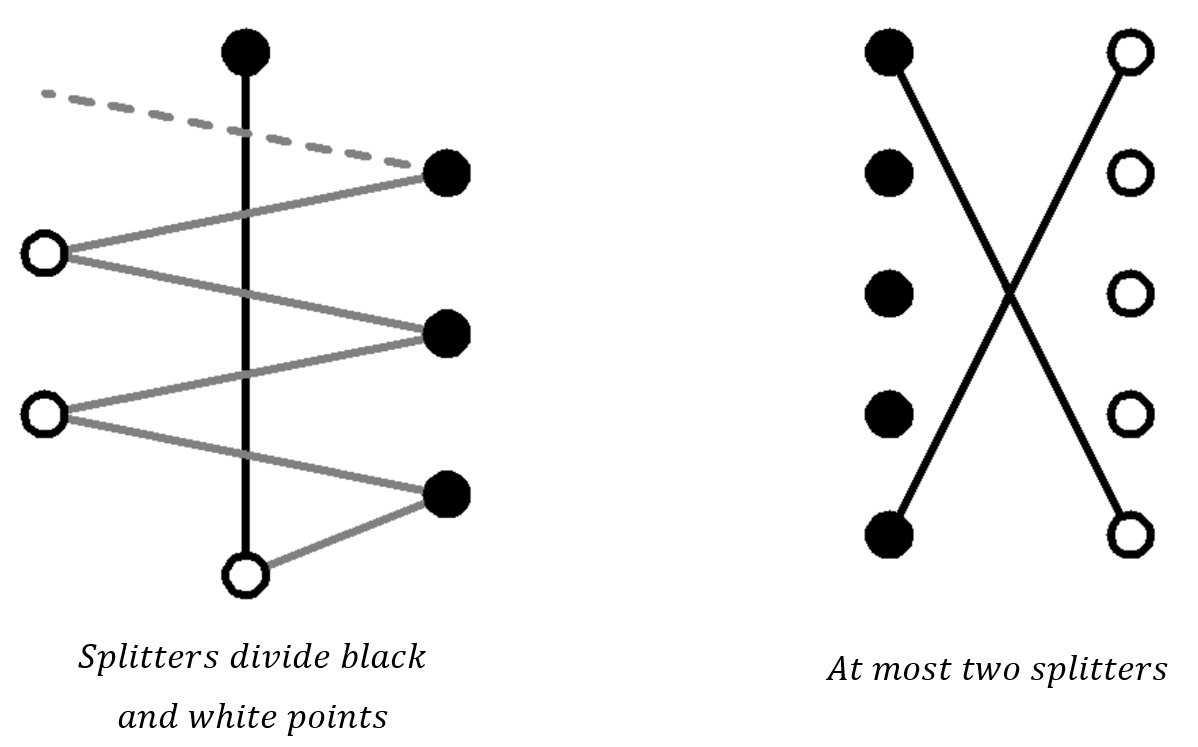}
\caption{Splitters divide bipartite graphs into colours, resulting in at most two splitters.}
\label{even_split}
\end{figure}

In the case that $N$ is even, we can no longer have all line-segments being splitters. This is because, for even $N$, a $N$-cycle graph is bipartite \cite{Wilson}. That is corners can be coloured black or white such that no line-segment connects two corners of the same colour. Consider a single splitter. The corners on either side of the line-segment must be different colours. Adding the subsequent line-segments results in all the black corners being on one side of the splitter, and all the white corners being on the other side. This means there can be at most two splitters in a $N$-cycle graph, for even $N$, as shown in Figure \ref{even_split}.

Therefore the best arrangement we can achieve, for $N$ even, is two splitters and $N-2$ one-off splitters. This is the arrangement we see for $N = 4$. 

\begin{figure}[H]
\centering
\includegraphics[width=0.8\textwidth]{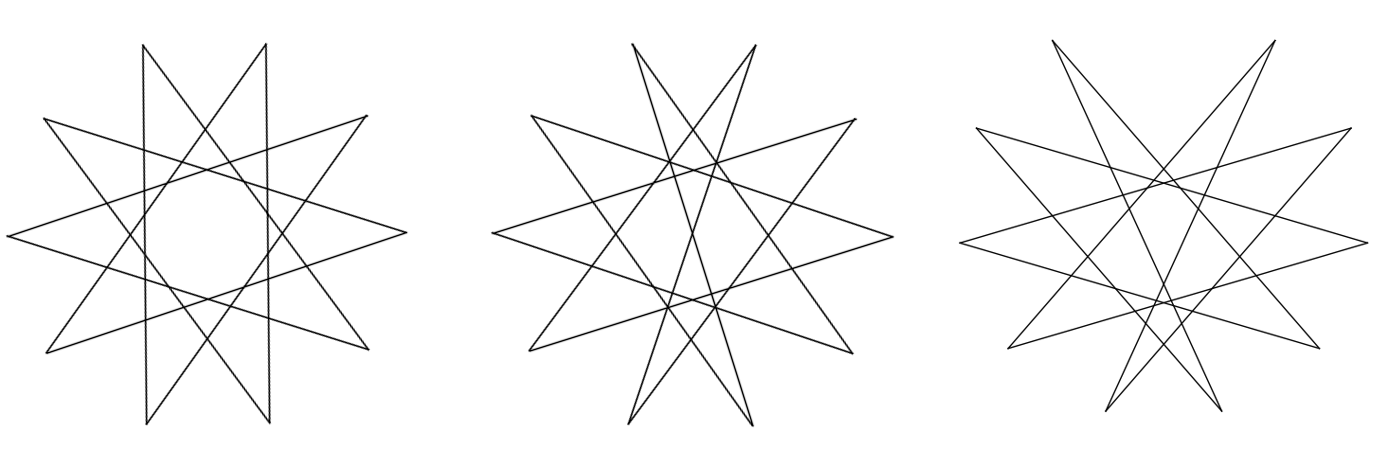}
\caption{Optimal arrangement for $N = 10$.}
\label{10}
\end{figure}

\begin{figure}[H]
\centering
\includegraphics[width=1\textwidth]{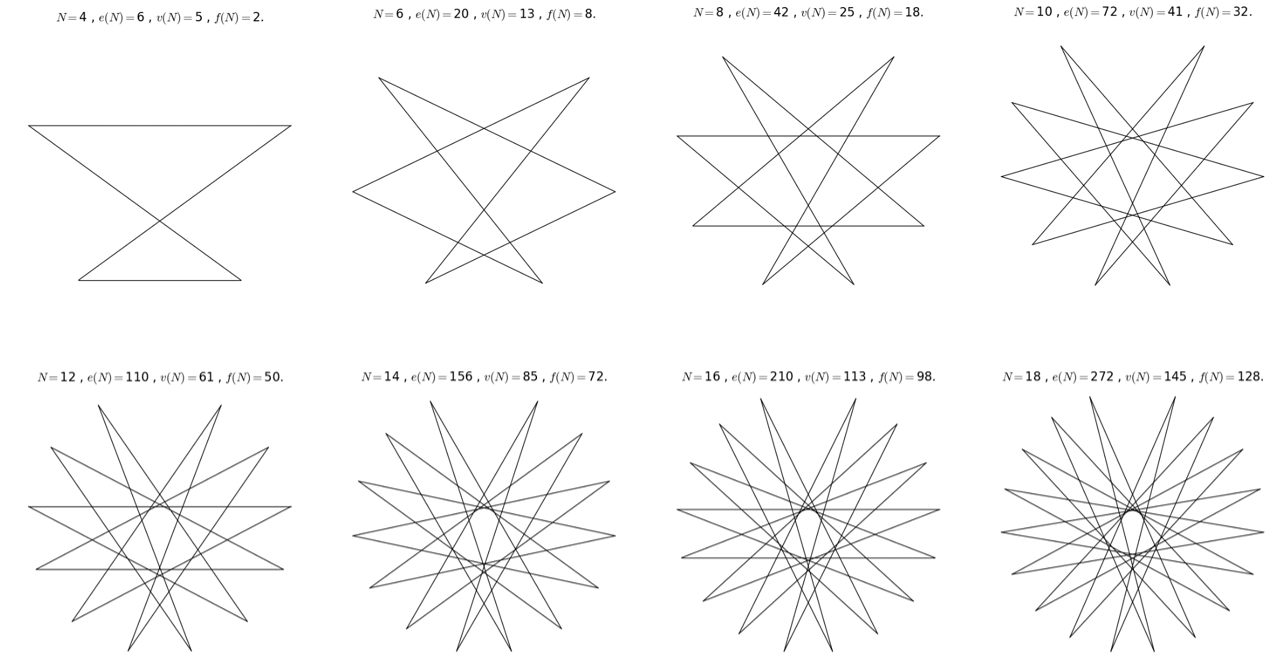}
\caption{Optimal arrangements for even $N$.}
\label{even}
\end{figure}

For even $N$, we can achieve this arrangement as follows. For all corners $1$ to $N$ connect corner $c$ to $\big(c + \frac{N-2}{2}\big)$ mod $N$. This results in $\frac{N}{2}$ pairs of parallel line-segments. Replace one of these pairs with two line-segments that cross each other, for example line-segments $\big(1 , \frac{N}{2}\big)$ and  $\big(\frac{N}{2} +1 , N\big)$ become $\big(1 ,\frac{N}{2} +1\big)$ and  $\big(\frac{N}{2}, N\big)$. Finally to eliminate most triple intersections arrange the corners on the vertices of a regular $N+1$-gon, with the blank space between the two splitters. This process is illustrated with the example $N=10$ in Figure \ref{10}. Note that for larger $N$ this does not eliminate all triple intersections, but perturbing the corners into a slightly irregular arrangement eliminates these.  
This results in the arrangements shown in Figure \ref{even}. 
We now use Euler's formula to work out the number of regions. We know that each of the $N-2$ line-segments contains $N-2$ vertices and two line-segments with $N-1$ vertices. Accounting for the fact that each vertex belongs to exactly two line-segments gives,
\[
v(N) = \frac{(N-2)(N-2) + 2(N-1)}{2}.
\]
As for edges, we have $N-2$ line-segments that contain $N-3$ edges and two line-segments with $N-2$ edges. This gives,
\[
e(N) = (N-2)(N-3) + 2(N-2).
\]
We then combine these using Euler's formula to give us 
\begin{align*}
  f(N) &= e(N) - v(N) + 1 \\
  &=  \frac{1}{2} N^2 - 2 N + 2.
\end{align*}
This concludes the proof of (\ref{eventhm}).


\begin{thebibliography}{9}

\bibitem{Diestel}
Reinhard Diestel (2017) \emph{Graph Theory}, Springer Berlin Heidelberg.

\bibitem{Wilson}
Robin J. Wilson (2010) \emph{Introduction to Graph Theory},  Pearson Education Limited.


\end{thebibliography}
\end{document}